\documentclass[preprint,12pt]{elsarticle}

\usepackage{amssymb}
\usepackage{amsmath}

\newtheorem{theorem}{Theorem}[section]

\newtheorem{corollary}{Corolarry}[section]

\numberwithin{equation}{section}

\numberwithin{equation}{section}
\journal{Journal of Computational and Applied Mathematics}

\begin{document}

\begin{frontmatter}



\title{On sufficient conditions in the classical problem of the calculus of variations}


\author[label1]{Misir \,J.~Mardanov} 

\affiliation[label1]{organization={Institute of Mathematics and Mechanics of Ministry of Science and Education},
            city={Baku},
            country={Azerbaijan}}
\author[label1,label2]{Telman.\,K.~Melikov} 

\affiliation[label2]{organization={Institute of Control Systems of Ministry of Science and Education},
            city={Baku},
            country={Azerbaijan}}

\author[label1,label3]{Samin\,T.~Malik} 

\affiliation[label3]{organization={ADA University},
            city={Baku},
            country={Azerbaijan}}

\begin{abstract}

This article is devoted to obtain new sufficient conditions for an extremum in problems of classical calculus of variations. 
The concept of a set of integrands is introduced.  Using this concept, first and second order sufficient conditions for a weak and strong local minimum, as well as an absolute minimum were obtained. Also, this concept, in particular, allows us to define a class of variation problems for which the necessary Weierstrass condition is also sufficient condition. It is shown that the sufficient conditions for a minimum obtained here have new areas of application compared to the known sufficient conditions of the classical calculus of variations. The effectiveness of the obtained results is illustrated by examples.

\end{abstract}


\end{frontmatter}



\section{ Introduction and problem statement }
\label{sec1}

Let us consider a classical problem of the calculus of variations
\begin{equation}
S(L)(x(\cdot))=\int_{t_{o}}^{t_{1}}L(t,x(t),\dot{x}(t))dt\rightarrow \mathop{min}_{x(\cdot)},  
\end{equation}

\begin{equation}
x(t_{o})=x_{0},\,x(t_{1})=x_{1},\quad \quad x_{0},\,x_{1}\in R^{n},
\end{equation}
where $R^{n}$ is $n$- dimensional Euclidean space,
 $x_{0},\,x_{1},\,t_{0}$ and $t_{1}$ are given points, and  
 $L(t,x,
\dot{x}): (a,\,b)\times R^{n}\times R^{n}\rightarrow R:=(-\infty,\,+\infty)$ is a given function, where $(a,\,b)-$ is some interval,
 and $[t_{0},\,t_{1}]\subset(a,\,b)$.
The function $L(\cdot)$ is called an integrand and assumed to be continuously differentiable with respect to the set of variables.  
The desired function $x(\cdot):I\rightarrow R^{n}$ is a piecewise-smooth vector function, where  $I:=[t_{0},\,t_{1}]$. It signifies that the function is continuous, and its derivative is continuous everywhere on $I=[t_{0},\,t_{1}]$ except for a finite number of points $\tau_{i}\in (t_{0},\,t_{1}),\,i=\overline{1,m}$. Moreover, at the point $\tau_{i}$ the derivative $\dot{x}(\cdot)$ has discontinuities of the first type (at the points $t_{0}$ and  $t_{1}$ the values of the derivative $\dot{x}(\cdot)$ are understood as the finite derivative on the right and left, respectively). In this case, as a rule, the points  $\tau_{i},\,\,i=\overline{1,m}$ are called the corner points of the function $\dot{x}(\cdot)$. The set of all piecewise smooth vector functions $x(\cdot):I\rightarrow R^{n}$ is denoted by $KC^{1}(I,\,R^{n})$.
 
Functions $x(\cdot)\in KC^{1}(I,\,R^{n})$, satisfying the boundary conditions (1.2) will be called admissible. The set of admissible functions will be denoted by $D$.   

Recall  (see, for example \cite[p. 107]{8}) some concepts that are more characteristic of the classical calculus of variations. An admissible function $\bar{x}(\cdot)$ is called a strong (weak) local minimum in the problem (1.1), (1.2), if there exists a number $\bar{\delta}>0\,(\tilde{\delta}>0)$ such that the inequality $S(L)(x(\cdot))\geq S(L)(\bar{x}(\cdot))$ is satisfied for all admissible functions $x(\cdot)$, for which

$$||x(\cdot)-\bar{x}(\cdot)||_{C(I, \, R^{n})}\leqslant\bar{ \delta}\,(max \lbrace ||x(\cdot)-\bar{x}(\cdot)||_{C(I, \, R^{n})},\,\,||\dot{x}(\cdot)-\dot{\bar{x}}(\cdot)||_{L_{\infty}(I, \, R^{n})}\rbrace\leqslant\tilde{\delta}).$$

It is clear that, if the admissible function $\bar{x}(\cdot)\in KC^{1}(I, \, R^{n}) $ provides a strong extremum, then it also provides a weak extremum. 

Therefore, for such functions, a necessary condition for a weak local extremum is a necessary condition for a  strong local one, and a sufficient condition for a strong local extremum is a sufficient condition for a weak local one.

It should be noted that for all  admissible functions $ x(\cdot)$ if the inequality  $S(L)(x(\cdot))\geq S(L)(\bar{x}(\cdot))$  is satisfied, then the admissible function $ \bar{x}(\cdot)$ is called the absolute minimum in problem (1.1), (1.2) and the notation $ \bar{x}(\cdot)\in \lbrace abs min \rbrace$ is adopted.

Also, let us recall some known necessary conditions for a minimum for the considered problem (see, for example, \cite{8, 2, 6, 12}).  Let the admissible function $ \bar{x}(\cdot)$ be a weak local minimum for problem (1.1), (1.2). In addition, let $\lbrace \tau \rbrace $ be the set of the corner points $\bar{x}(\cdot)$. Then:

(i) the function $\bar{x}(\cdot)$ satisfies the following equations, i.e. the equalities are true:

(a) \begin{equation}
\delta S (L)(\Delta x(\cdot);\, \bar{x}(\cdot))=\int_{t_{0}}^{t_{1}}[\bar{L}_{x}^T(t)\Delta x(t)+\bar{L}_{\dot{x}}^T(t)\Delta \dot{x}(t)]dt=0,
\end{equation}
for all admissible functions $x(\cdot)$, where $\Delta x(t)=x(t)-\bar{x}(t), \,t\in I$;

(b) there exists a constant vector $ c\in R^{n}$ such that

\begin{equation}
\bar{L}_{\dot{x}}(t)- \int_{t_{0}}^{t} \bar{L}_{x}(\nu)d\nu=c, \,\, t\in I\backslash \lbrace \tau \rbrace,
\end{equation}
here $\delta S (L) (\Delta x(\cdot);\, \bar{x}(\cdot))-$ is the first variation of the functional (1.1),  the equation (1.4) is the Euler equation in integral form, where $$\bar{L}_{x}(\mu)=L_{x}(\mu, \bar{x}(\mu),\dot{\bar{x}}(\mu)),\,\bar{L}_{\dot{x}}(\mu)=L_{\dot{x}}(\mu, \bar{x}(\mu),\dot{\bar{x}}(\mu)), \,\, \mu\in\lbrace t,\,\nu\rbrace,$$ is the derivative of the integrand $L(t,x,y)$;

(ii) for each $ \tau \in \lbrace  \tau \rbrace $ along $\bar{x}(\cdot)$ the Weiestrass-Erdman conditions are satisfied, i.e. the equalities are true
\begin{equation}
\bar{L}_{\dot{x}} (\tau-)=\bar{L}_{\dot{x}}  (\tau+),\quad \quad \bar{L} (\tau-)-\dot{\bar{x}}^{T}(\tau-) \bar{L}_{\dot{x}} (\tau-)=\bar{L}(\tau+)-\dot{\bar{x}}^{T}(\tau+) \bar{L} _{\dot{x}} (\tau+),
\end{equation}
here and in what follows, for example, the symbols $\bar{L }(\tau+)\,\,(\bar{L} (\tau-)) $ mean the limit from the right (left) of the function $\bar{L}(t)$ at the point $t=\tau$;

(iii) By virtue of  \cite{14}, \cite{15} there exists a number $\delta>0$, for which along function $\bar{x}(\dot)$ the inequalities hold

\begin{align*}
E(\bar{L})(t,\, \xi)\geq0,\quad E(\bar{L}) (\tau\pm,\, \xi)\geq0,
\end{align*}
\begin{align}
\forall t\in I\setminus\lbrace\tau\rbrace,\quad \forall\tau\in\lbrace\tau\rbrace, \quad\forall \xi\in B_{\delta}(0).
\end{align}

Here $B_{\delta}(0)-$ is a closed ball of radius $\delta$ with the center $0\in R^{n}$, further, 

\begin{equation*}
E(\bar{L})(t,\, \xi):= E(L)(t, \bar{x}(t), \dot{\bar{x}}(t),\dot{\bar{x}}(t)+ \xi)=
\end{equation*}
\begin{equation}
=L(t, \bar{x}(t), \dot{\bar{x}}(t)+ \xi)-L(t, \bar{x}(t), \dot{\bar{x}}(t))-\bar{L}^{T}_{\dot{x}}(t)\xi,
\end{equation}
where $ E(L)(t,x,y,z)-$ the Weierstrass function for the problem (1.1), (1.2), which is defined as 
\begin{equation}
E(L)(t,x,y,z)=L(t,x,z)-L(t,x,y)-L_{y}^{T}(t,x,y)(z-y).
\end{equation}

The following statement is also known  from \cite{2}, \cite{6}.
If the admissible function $\bar{x}(\cdot)$ is a strong local minimum, then the Weierstrass condition is satisfied, i.e. the inequalities are valid 

\begin{align*}
E(\bar{L})(t,\, \xi)\geq0,\quad \forall t \in I\setminus\lbrace\tau\rbrace, \quad\forall \xi\in R^{n}, 
\end{align*}
\begin{align}
E(\bar{L})(\tau\pm,\, \xi)\geq0,\quad \forall \xi\in R^{n}.
\end{align}

We also note that from (1.6) the Legendre conditions follow as a consequence. We formulate this condition \cite{2}, \cite{6}. If the admissible function $\bar{x}(\cdot)$ is a weak local minimum in problem (1.1), (1.2) in addition, integrand $L(t,x,\dot{x})$  is twice differentiable with respect to the variable $\dot{x}$, then for all $\xi\in R^{n}$ the following inequalities hold

\begin{align*}
\xi^{T}\bar{L}_{\dot{x}\dot{x}}(t)\xi\geq0,\quad \forall t \in I\setminus\lbrace\tau\rbrace,  
\end{align*}
\begin{align}
\xi^{T}\bar{L}_{\dot{x}\dot{x}}(\tau\pm)\xi\geq0,\quad \forall \tau\in \lbrace\tau\rbrace.
\end{align}

Each function $\tilde{x}(t)\in KC^{1}(I,\,R^{n})$, satisfying the integral equation (1.4) or, i.e., the equation (1.3), is usually called an extremal of the problem (1.1), (1.2). If the extremal is an element of the set $D$, then it is called an admissible extremal. 

The main goal of this paper is to study the minimum of admissible extremals of the problem (1.1), (1.2). 

Let $\bar{x}(\cdot)\in D-$ be some extremal of the problem (1.1), (1.2). If along it at least one of the necessary conditions (1.5), (1.6) and (1.9) is not satisfied, then it is obvious that the extremal $\bar{x}(t)\in D$ is not even a weak local minimum in the problem (1.1), (1.2). So, for this case, the study of the problem (1.1), (1.2) is complete. Therefore, it is assumed below that along the extremal $\bar{x}(\cdot)\in D $ the conditions (1.5), (1.6) and (1.9) are satisfied.

Consequently, it is expedient to check the fulfillment of the necessary conditions (1.5), (1.6), and (1.9) along the extremal $\bar{x}(\cdot)\in D $ under consideration.

It should be noted that obtaining sufficient conditions for a minimum in the theory of variational calculus is one of the main directions.

The question of sufficient conditions in the calculus of variations was first \linebreak systematically studied by Karl Jacobi \cite{9}. Introducing the concept of a conjugate point, he obtained a necessary and sufficient condition for a weak local extremum.

The study of a strong local extremum in the theory of variational calculus is associated with the name of Karl Weierstrass. He established a necessary and close to it a sufficient condition for a strong extremum.
Weierstrass's research in the calculus of variations was then developed by Adolf Mayer \cite{18}, Oskar Bolza\cite{3} and others.

The true essence of the theory of strong extremum in the calculus of variations was revealed by David Hilbert and Henri Poincaré.

In the work \cite[pp. 299-301]{8}, developing Hilbert's idea \cite{7}, Carathéodory \cite{4} and Young \cite{19} introduced 
the concept of a local $ K -$ function for problems in the classical calculus of variations. Global $ K -$ functions were 
considered by Krotov (\cite{10}, \cite{11}) and Risone Maicho (\cite{16}, \cite{17}). Further, a meaningful study of sufficient 
conditions was carried out in the works of Levitin, Milyutin, and Osmolovsky
 \cite{13}. The idea of the work \cite{13} is  very close to the approach
 carried out in the work \cite[pp. 305-314]{8}.

Recall (see \cite{2}, \cite{6}) that in classical variational problems, a number of known sufficient conditions for a minimum, in general, remain powerless when studying an extremal if at some point $ (\theta,\, \eta)\in I\times R^{n}\setminus\lbrace0\rbrace$ at least one of the Weierstrass and Legendre conditions degenerates, i.e. at the point $ (\theta,\, \eta)\in I\times R^{n}\setminus\lbrace0\rbrace$  at least one inequality from (1.8) and (1.9) is satisfied as an equality.

A comparison of the above results on the sufficiency of minimum conditions for the problem (1.1), (1.2) shows that they have different shortcomings and also different areas of application (see example 5.1). 

Therefore, we can say: 1) the theory of sufficient conditions for problems of the calculus of variations is far from complete;

2) obtaining new sufficient conditions for the problem (1.1), (1.2) is still relevant today.

In this article, following the property of invariance of the Hilbert integral (\cite [pp. 326-335 ]{8}), an approach is proposed to obtain new sufficient conditions for the minimum of the extremal of the problem (1.1), (1.2).

In recent years, for example, in the works \cite{20}, \cite {21}: (1) various sufficient conditions of the Jacobi type have been obtained; (2) verification of the fulfillment of the Jacobi condition is reduced to investigating the existence of a continuous solution to the Riccati matrix equation (this study is called the Riccati method);
(3) it has been shown that the Riccati method is effective for applying a modern computational approach to solving, for example, problem (1.1), (1.2).

As can be seen, the statements of Theorems 4.1 and 4.2 are expressed in the form of differential inequalities (4.1) and (4.2), respectively. This means that the application of Theorems 4.1 and 4.2 leads to the determination of solutions to the differential inequalities (4.1) and (4.2), respectively.

The works  (\cite {27}, \cite {26}), for example, are devoted to the study of differential inequalities. Following the above, it can be stated that the study of differential inequalities remains relevant today as well.

It should be noted that some problems of the calculus of variations and the theory of optimal control are studied in more detail in the work of \cite{28}.

The structure of the present paper is outlined by the following scheme. In second section, The concept of a set of integrands for problem (1.1), (1.2) is introduced and some properties of this concept are studied.
In third section are obtained some formulas for the increment of a functional in terms of a set of integrands.
In the third section, various formulas are obtained for the increment of the functional in terms of a set of integrands.

In the fourth section, sufficient conditions for a minimum for the extremal of problem (1.1), (1.2) are obtained.

The last section, using specific examples, demonstrates the effectiveness of the results obtained in the fourth section.

\section {The concept of a set of integrands for the problem (1.1), (1.2).} 

We introduce following sets of the form for the problem (1.1), (1.2) and denote it by $ Z (L)$. Namely:
\begin{equation}
Z(L)=\lbrace L^{*}(t,x,\dot{x})\in C^{1}(I\times R^{n}\times R^{n}, R): S(L^{*})(x(t))=S(L)(x(t)),\,\, \forall x (\cdot)\in D \rbrace
\end{equation}
where $D-$ set of admissible functions. 

\textbf{Definition 2.1.} The set defined by (2.1), is called the set of integrands for the problem (1.1), (1.2).

This means, $ L(\cdot)\in Z( L(\cdot))$. Let $\bar{x}(t)\in D $. Then it is clear that the set of the form 

 $$\lbrace L^{*}(t,x(t),\dot{x}(t))= L(t,x(t),\dot{x}(t))$$
 $$+\frac{d}{dt}[q(t)(x(t)-\bar{x}(t)],\,\, \forall x(t) \in D, \, \forall q(t)\in KC^{1}(I,R^n) \rbrace, $$
 is a subset of the set $Z(L)$. This means, $Z(L)\backslash\lbrace L(\cdot)\rbrace\neq\varnothing$.

We especially note that the concept of a set of integrands (see Definition 2.1), without requiring an assumption about the field of extremals (see \cite{2}, \cite{6}), allows us to obtain a strong and weak local and also absolute minimum in the problem (1.1), (1.2).

Let us consider the sets of problem (1.1), (1.2) obtained by replacing the integrand $L^{*}(\cdot)$ of each element by the set  $Z(L)$. It is clear that, by Definition 2.1, for such problem the set of extremal coincides and an inequality of the form is satisfied
\begin{equation*}
\int^{t_{1}}_{t_{0}}L^{*}(t,x(t),\dot{x}(t))dt\geq\int_{t_{0}}^{t_{1}}L(t,\, \bar{x}(t)\,\dot{\bar{x}}(t))dt,
\end{equation*}

\begin{equation}
\forall x(t) \in D, \,\,\forall L^{*}(t,x,\dot{x})\in Z(L).
\end{equation}
where the function  $\bar{x}(\cdot)\in D$ is minimum in the problem (1.1), (1.2). From here, by virtue of Definition 2.1, we obtain the statement:
if the inequality (2.2) is satisfied, it is a necessary and sufficient condition for the minimum of the extremal $\bar{ x}(\cdot)\in D$ in the problem (1.1), (1.2).

It is interesting to note that inequality (2.2) does not, in general, imply an inequality of the form

\begin{equation}
L^{*}(t,\, x(t),\,\dot{x}(t))\geqslant L(t,\, \bar{x}(t),\,\dot{\bar{x}}(t)),
\forall x(\cdot)\in D,\,\text{a.e.}\, \,\, t\in I,  
\end{equation}
\begin{equation*}
  \forall L^{*}(t,x,\dot{x})\in Z(L),  
\end{equation*}
where $\bar{x}(\cdot)\in D-$ is an admissible extremal of the problem (1.1), (1.2).

However, there may exist an integrand $L^{*}(\cdot)\in Z(L)$ such that for the extremal $\bar{x}(\cdot)\in D$ the inequality holds

\begin{equation}
L^{*}(t,\, x(t),\,\dot{x}(t))\geqslant L(t,\, \bar{x}(t),\,\dot{\bar{x}}(t)),
\forall x(\cdot)\in D,\,\text{a.e.} \,\,\, t\in I.
\end{equation}

To confirm the last stated statements (see (2.3), (2.4)) it is sufficient to consider the following example.

\textbf{Example 2.1}.
\begin{equation}
\int_{0}^{1}{\dot{x}}^{2}dt \longrightarrow min,\,\, x(0)=1, x(1)=0.
\end{equation}

Here $L(t,x,\dot{x})=\dot{x}^{2}$ and $\bar{x}(t)=1-t,\,\,t\in [0, 1]-$ is an admissible extremal in the problem (2.5). First we show that inequality  (2.3) does not hold for the problem (2.5), when $\bar{x}(t)=1-t,\,\, t\in[0, 1]$. It is clear that inequality (2.3) for the problem (2.5) takes the form:  
$$(\dot{\bar{x}}(t)+\Delta\dot{x}(t))^{2}\geq\dot{\bar{x}}^{2}(t),\,\, \forall\Delta x(t)=x(t)-\bar{x}(t),\,x(t)\in D,$$
$$\Delta x(0)= \Delta x(1)=0,\, t\in [0, 1].$$

From here we have that inequality (2.3) is not satisfied, since the inequality $\Delta \dot{x}(t)[\Delta \dot{x}(t)-2]\geq 0$ is not satisfied for $\Delta \dot{x}(t)\in (0,\,2), \, t\in [0,1]$.

Now we will show that inequality (2.4) is satisfied.
 Let's choose $L^{*}(\cdot)$ in the following form:
$$L^{*}(t,x(t),\dot{x}(t))=L(t,x(t),\dot{x}(t))+2\frac{d}{dt}(x(t)-\bar{x}(t)),\,x(\cdot)\in D$$ 
where $L(\cdot)=\dot{x}^{2},\,\,\bar{x}(t)=1-t,\,\,t\in[0,\,1]$.
\\

It is obvious that  $L^{*}(\cdot)\in Z (L)$. In this case, the inequality (2.4) takes the form:
$$\dot{x}^{2}(t)+2(\dot{x}(t)-\dot{\bar{x}}(t))\geq \dot{\bar{x}}(t)^{2},\,\,\forall x(t)\in D,\,t\in[0,\,1].$$

Considering $\dot{\bar{x}}(t)=-1,\,\,t\in[0,\,1]$ we will get
$$(\dot{x}(t)+1)^{2}\geq0,\,\, \forall x(t)\in D,\,t\in[0,\,1].$$

Hence, for the problem (2.5) the inequality (2.4) is satisfied.

It is clear that by Definition 2.1 the inequality
$$\Delta S(L)(\bar{x}(\cdot))=S(L)(x(\cdot))- S(L)(\bar{x}(\cdot))\geq0,\,\,\forall x(\cdot)\in D.$$
                                                    
From the last inequality it follows that the admissible extremal  $\bar{x}(t)=1-t,\,\, t\in [0,\,1] $ is an absolute minimum in problem (2.5).

Therefore, based on example 2.1 and the above reasoning, we come to the conclusion: the usage of the concept of a set of integrands allows us to obtain sufficient conditions for a minimum in the problem (1.1), (1.2). 

\section {Some formulas for the increment of a functional in terms of a set of integrands.}  

Let the function $\bar{x}(\cdot)\in D$ be an extremal in the problem (1.1), (1,2) and $x(t)=\bar{x}(t)+\Delta x(t)\in D$ -- be an arbitrary function, where $\Delta x(t_{0})=\Delta x(t_{1})=0$.

Let us consider the set of integrands relative to $\bar{x}(\cdot)$ in the form

\begin{align*}
Z(L)(\bar{x}(\cdot))=\lbrace L^{*}(t,x,\dot{x})\in Z (L) : L^{*}(t, x(t),\dot{x}(t))=L(t, x(t),\dot{x}(t))+ \\
+\frac{1}{2} \frac{d}{dt}[\Delta x^{T}(t)Q(t)\Delta x(t)],\,\, \text{where} \,\, x(t)=\bar{x}(t)+\Delta x(t),\,\, t\in I,
\end{align*}
\begin{align}
\quad\quad \Delta x(t_{0})=\Delta x(t_{1})=0,\,\, \Delta x(t)\in KC^{1}(I,\,R^{n})\rbrace.
\end{align}

Here and below $ n\times n $ matrix function $Q(t),\,t\in [t_{0},\,t_{1})$ is an unknown function that satisfies the following assumptions
\begin{align*} 
Q^{T}(t)=Q(t),\,t\in [t_{0},\,t_{1}),\, \,Q(\cdot)\in KC^{1}([t_{0},\,t_{1}),\,R^{n\times n}),
\end{align*}
\begin{equation}
\lim_{t\rightarrow t_{1}}\Delta x^{T}(t)Q(t)\Delta x(t)=0,\,\,\,\lim_{t\rightarrow t_{1}}\Delta x^{T}(t)Q(t)\Delta \dot{x}(t)\in (-\infty, +\infty), 
\end{equation}

$\mathop{\lim}\limits_{t\rightarrow t_{1}}\Delta x^{T}(t)\dot{Q}(t)\Delta x(t)\in (-\infty, +\infty)$, i.e. there exists finite limits.

Obviously,
\begin{equation}
Z(L)(\bar{x}(\cdot))\subset Z(L).
\end{equation}

Let $L^{*}(t,x, \dot{x})\in Z(L)(\bar{x}(\cdot))$. Then, by definition 2.1, taking into account (3.1)-(3.3) and $x(t)=\bar{x}(t)+\Delta x(t)$, where $\Delta x(t_{0})=\Delta x(t_{1})=0$ for the increment $\Delta S(L)(\bar{x}(\cdot))$ the following equality is valid:
\begin{align}
\begin{split}
\Delta S(L)(\bar{x}(\cdot))=\int_{t_{0}}^{t_{1}}L(t,x(t),\dot{x}(t))dt-\int_{t_{0}}^{t_{1}}L(t,\bar{x}(t),\dot{\bar{x}}(t))dt=\\
=\int_{t_{0}}^{t_{1}}L^{*}(t,x(t),\dot{x}(t))dt-\int_{t_{0}}^{t_{1}}L(t,\bar{x}(t),\dot{\bar{x}}(t))dt=\\
=\int_{t_{0}}^{t_{1}}[L(t,\bar{x}(t)+\Delta x(t),\dot{\bar{x}}(t)+\Delta \dot{x}(t))+\Delta \dot{x}^{T}(t)Q(t)\Delta x(t)\\
+ \frac{1}{2}\Delta x^{T}(t)\dot{Q}(t)\Delta x(t)]dt -\int_{t_{0}}^{t_{1}}L(t,\bar{x}(t),\dot{\bar{x}}(t))dt.
\end{split}
\end{align}

Let us obtain from (3.4) the following increment formulas for the problem (1.1), (1.2):

(i) let the integrand $L(\cdot)$, its partial derivatives $ L_{x}(\cdot)$ and $ L_{\dot{x}}(\cdot)$ be continuous with respect to the set of variables. Then, applying the Taylor formula with the Lagrange remainder term from (3.4), we obtain

\begin{align}
\begin{split}
\Delta S(L)(\bar{x}(\cdot))=\int_{t_{0}}^{t_{1}}\{L(t,\bar{x}(t), \dot{\bar{x}}(t)+\Delta \dot{x}(t))-L(t,\bar{x}(t), \dot{\bar{x}}(t))+\\
+[L_{x}(t,\bar{x}(t)+\theta\Delta x,\dot{\bar{x}}(t)+\Delta \dot{x}(t))+ \Delta \dot{x}(t) Q(t)]\Delta x(t)+\\
+\frac{1}{2}\Delta x^{T}(t)\dot{Q}(t)\Delta x(t)\}dt,
\end{split}
\end{align}
where $\theta \in (0,\,1),\,\, \bar{x}(\cdot)\in D$ is extremal in the problem (1.1), (1.2).

We use (1.7), i.e. the Weierstrass function, and take into account equalities (1.3). Then the equalities (3.5) take the following form

\begin{align}\label{3.6}
\begin{split}
\Delta S(L)(\bar{x}(\cdot))=\int_{t_{0}}^{t_{1}}\{E(L)(t,\bar{x}(t), \dot{\bar{x}}(t),  \dot{\bar{x}}(t)+\Delta \dot{\bar{x}}(t))+[\Delta \dot{x}^{T}(t)Q(t)\\
+L_{x}(t,\bar{x}(t)+\theta\Delta x(t), \dot{\bar{x}}(t)+\Delta \dot{x}(t))- L_{x}(t,\bar{x}(t), \dot{\bar{x}}(t))]\Delta x(t)\\
+\frac{1}{2}\Delta x^{T}\dot{Q}(t)\Delta x(t)\}dt;
\end{split}
\end{align}

(ii) let the function  $L(t,x,\dot{x}) $ and its partial derivatives $ L_{x}(\cdot),L_{\dot{x}}(\cdot)$ and  $L_{xx}(\cdot),$  be continuous in the set of variables on $ I\times R^{n}\times R^{n}$. Then applying the Taylor formula with the Lagrange remainder term from (3.4) we have
  
\begin{align}
\begin{split}
\Delta S(L)(\bar{x}(\cdot))=\int_{t_{0}}^{t_{1}}[L(t,\bar{x}(t), \dot{\bar{x}}(t)+\Delta \dot{x}(t))-L(t,\bar{x}(t), \dot{\bar{x}}(t))+\\
L_{x}(t,\bar{x}(t), \dot{\bar{x}}(t)+\Delta \dot{x}(t))\Delta x(t)\\
 +\frac{1}{2}\Delta x^{T}(t)L_{xx}(t,\bar{x}(t)+\theta \Delta x(t), \dot{\bar{x}}(t)+\Delta \dot{x}(t))\Delta x(t)\\
 +\frac{1}{2}\Delta x^{T}(t)\dot{Q}(t)\Delta x(t)+\Delta \dot{x}^{T}(t)Q(t)\Delta x(t)]dt, 
 \end{split}
\end{align}
where $\theta \in(0,\,1)$ and $\bar{x}(t)\in D-$ is the extremal in the problem (1.1), (1.2). 

Since the function $\bar{x}(t)\in D$ is an extremal in the problem (1.1), (1.2), then by virtue of the equality (1.3) and the notation (1.7) (i.e. by virtue of the Weierstrass function) equality (3.7) takes the form

\begin{align}\label{3.8}
\begin{split}
\Delta S(L)(\bar{x}(\cdot))=\int_{t_{0}}^{t_{1}}\{E(L)(t,\bar{x}(t), \dot{\bar{x}}(t), \dot{\bar{x}}(t)+\Delta \dot{x}(t))\\
+[\Delta \dot{x}(t)Q(t)+L_{x}(t,\bar{x}(t), \dot{\bar{x}}(t)+\Delta \dot{x}(t))
- L_{x}(t,\bar{x}(t), \dot{\bar{x}}(t)]\Delta x(t)\\+\frac{1}{2}\Delta x^{T}(t)[\dot{Q}(t)+L_{xx}(t,\bar{x}(t)+\theta \Delta x(t),\dot{\bar{x}}(t)+\Delta\dot{x}(t))] \Delta x(t)\}dt, 
\end{split}
\end{align}
where $E(L)(\cdot)-$ Weierstrass function (see (1.7)), and  $Q(t)-$ $n\times n$ matrix function that satisfies assumptions (3.2).

\section {Sufficient conditions for a minimum for the extremal of problem (1.1), (1.2).}

Using the increment formulas (3.6), taking into account the concept of a set of integrands (see Definition 2.1) and the approach given in solving Example 2.1, the following theorem is easily proven.

\begin{theorem}
Let $\bar{x}(t)\in D$ be an extremal of the problem (1.1), (1.2), and
let the integrand $L(\cdot)$, its partial derivatives $ L_{x}(\cdot)$ and $ L_{\dot{x}}(\cdot)$ be continuous in the set of variables on $I\times R^{n}\times R^{n}$ and there exists an $n\times n$ matrix function $Q(t),\,\,t\in [t_{0},\,t_{1})$, for which assumption (3.2) holds. If the inequality

\begin{align*}
E(L)(t,\bar{x}(t), \dot{\bar{x}}(t), \dot{\bar{x}}(t)+\xi)+[\xi^{T}Q(t)+L_{x}(t,\bar{x}(t)+\theta \eta,\, \dot{\bar{x}}(t)+\xi )-
\end{align*}
\begin{align}
-L_{x}(t,\bar{x}(t),\dot{\bar{x}}(t)]\eta+\frac{1}{2}\eta^{T}\dot{Q}(t)\eta]\geq0,\,\, t\in[t_{0},\,t_{1})\setminus\lbrace\tilde{t}\rbrace, \forall \theta \in (0,1)
\end{align}

\begin{align*}
\forall \xi \in R^{n}, \,\,\forall \eta \in R^{n}\,\, (\forall\xi\in R^{n},\,\eta \in B_{\delta}(0)\,\,(\forall\xi\in B_{\delta}(0),\,\forall \eta \in B_{\delta}(0))),
\end{align*}
 holds, then the extremal $\bar{x}(t)\in D $ is an absolute (strong local (weak local)) minimum in the problem (1.1), (1.2), where  $E(L)(\cdot)-$ Weierstrass function which is determined by (1.7),  $\lbrace\tilde{t}\rbrace\subset(t_{0},\,t_{1})-$ is a finite set, $B_{\delta}(0)-$ closed ball of radius $\delta$ with center $0\in R^{n}$.
\end{theorem}
Now, similar to Theorem 4.1, using the increment formula (3.8), the following statement is proved.

\begin{theorem} Let $\bar{x}(t)\in D$ be an extremal of the problem (1.1), (1.2) and the following assumptions are satisfied:

(i)	the integrand $L(\cdot):I\times R^{n}\times R^{n}\rightarrow R$ and its partial derivatives $L_x(\cdot),$ $L_{\dot{x}}(\cdot)$ and   $L_{xx}(\cdot)$  are continuous with respect to the set of variables;

(ii)	there exists a symmetric $n\times n$ matrix function $Q(t),\, t\in [t_{0},\,t_{1})$, for which (3.2) is satisfied;

(iii) the following inequality holds
\begin{align}
\begin{split}
E(L)(t,\bar{x}(t), \dot{\bar{x}}(t), \dot{\bar{x}}(t)+\xi)\\
+[\xi^{T}Q(t)+L_{x}(t,\bar{x}(t),\dot{\bar{x}}(t)+\xi)-L_{x}(t,\bar{x}(t),\dot{\bar{x}}(t)]\eta+\\
+\frac{1}{2}\eta^{T}[\dot{Q}(t)+L_{xx}(t,\bar{x}(t)+\eta\theta,\,\,\dot{\bar{x}}(t)+\xi)]\eta\geq 0,\\
\,\,\forall \theta \in (0,\,1),\,\,\forall t\in [t_{0},\,t_{1})
\setminus\lbrace\tilde{t}\rbrace,\,\,\,\,\,
\forall (\xi,\,\eta) \in R^{n}\times R^{n},\\ \,\,(\forall\xi\in R^{n},\, \forall \eta \in B_{\delta}(0))\,\,(\forall(\xi,\,\eta)\in B_{\delta}(0)\times B_{\delta}(0))),\\
\end{split}
\end{align}
where $\lbrace\tilde{t}\rbrace\subset (t_{0},\,t_{1})-$finite set, $B_{\delta}(0)-$ closed ball of radius $\delta$ with center $0\in R^{n}$.

Then the extremal $\bar{x}(\cdot)\in D$ is an absolute (strong local (weak local)) minimum in the problem (1.1), (1.2).
\end{theorem}

From Theorem 4.1 and Theorem 4.2, taking into account the increment formulas \eqref{3.6} and (3.8), the following statements directly follow, respectively.

\begin{theorem}

Let $\bar{x}(\cdot)\in D$ be an extremal of the problem (1.1), (1.2) and the following assumptions are fulfilled:

(i) the integrand $L(\cdot)$ has the form 
\begin{align*}
L(t,x,\dot{x})=L^{(0)}(t,\dot{x})+q^{T}(t)x, \,\,\, 
\end{align*}
where  the function $L^{(0)}(\cdot):I\times R^{n}\rightarrow R$  and its partial derivative $L^{(0)}_{\dot{x}}$, are continuous, and $q(t)\in C(I)$;

(ii) the inequality
\begin{align}
E(L)(t,\bar{x}(t), \dot{\bar{x}}(t), \dot{\bar{x}}(t)+\xi)\geq0,\,\, \forall t \in I \setminus\lbrace\tau\rbrace,\,\,\forall \xi \in R^{n}\,\,\,(\forall \xi \in B_{\delta}(0)) 
\end{align}
holds. Then the extremal $\bar{x}(\cdot)\in D$ is an absolute (weak local) minimum in the problem (1.1), (1.2).
\end{theorem}

\begin{theorem}
Let $\bar{x}(\cdot)\in D$ be an extremal of the problem (1.1), (1.2) and the following assumptions are satisfied:

(i)	the integrand $L(t,x,\dot{x})$ has the form: $L(t,x,\dot{x})=L^{(0)}(t,\dot{x})+L^{(1)}(t,\,x)$, where the functions $L^{(0)}(t,\dot{x})$,  $L^{(1)}(t,\,x)$ and their partial derivatives 
$L^{(0)}_{\dot{x}}(\cdot),$ $\,L^{(1)}_{x}(\cdot),\,L^{(1)}_{xx}(\cdot)$ are continuous in the set of variables;

(ii) the following inequalities hold
\begin{align}
E(L^{(0)})(t,\dot{\bar{x}}(t), \dot{\bar{x}}(t)+\xi)\geq0,
\end{align}

\begin{align}
\eta^{T}L^{(1)}_{xx}(t,\bar{x}(t)+\theta\eta)\eta\geq0,\,\,t \in I \setminus\lbrace\tau\rbrace,
\end{align}

$\forall\theta \in (0,1),\,\forall( \xi,\,\eta) \in R^{n}\times R^{n}\,\,\,(\forall (\xi,\,\eta) \in  R^{n}\times B_{\delta}(0))\,(\forall(\xi,\,\eta)\in B_{\delta}(0)\times B_{\delta}(0)))$, where $\lbrace\tau\rbrace$ is the set of points of the corner of the extremal $\bar{x}(\cdot)\in D$.

Then the extremal $\bar{x}(\cdot)\in D$ is an absolute (strong local (weak local)) minimum in the problem (1.1), (1.2).
\end{theorem}

The following statement follows directly from Theorem 4.4.

\begin{corollary} 

Let the integrand $L(\cdot)$ be linear with respect to the variable $x$, i.e. $L(t,x,\dot{x})=L^{(0)}(t,\dot{x})+q^{T}(t)x$. 

Then the fulfillment of the Weierstrass condition (4.3) for the extremal $\bar{x}(\cdot)\in D$ is a necessary and sufficient condition for an absolute (weak local) minimum in the problem (1.1), (1.2).
\end{corollary}

The proof of the corollary follows directly from Theorem 4.4, taking into account (1.6) and (1.8).

\section{Discussion and Examples}

As can be seen, the statements of Theorems 4.1 and 4.2 are expressed in
the form of differential inequalities (4.1) and (4.2), respectively. This means
that the application of Theorems 4.1 and 4.2 leads to the determination of
solutions to the differential inequalities (4.1) and (4.2), respectively.
Let us note that, for example, the works (\cite {27}, \cite {26})  are devoted to the study of differential inequalities.
Following the above, it can be stated that the study of differential
inequalities remains relevant today as well.
In some special cases, for example, for the case of the form $L(t,x, \dot{x})=<A(t)\dot{x},\dot{x}>+2<B(t)\dot{x},x>+<C(t)x,x>$, solutions
of the differential inequality (4.2) lead to solutions of the Riccati matrix
equation.

To demonstrate the effectiveness of the results obtained in Section 4, consider the following examples.
 
\textbf{Example 5.1.}

\begin{equation}
\int^{1}_{-1}(t^{2}\dot{x}^{2}+12x^{2})dt\rightarrow min,\,\,\, x(-1)=-1,\,\,\,x(1)=1,
\end{equation}
where $L(t,x,\dot{x})=t^{2}\dot{x}^{2}+12x^{2}$. 

It is clear that the Euler equation in problem (5.1) has the form:

\begin{equation}
t^{2}\ddot{x}+2t\dot{x}-12x=0,\,\,t\in[-1,\,1],
\end{equation}

As is known, this equation is an Euler type equation \cite[p. 110]{5}.

Using \cite[p. 111]{5}, we obtain that the general integral of the equation (5.2) has the form:

$$x(t)=c_{1}t^{3}+c_{2}\frac{1}{t^{4}},\,\,\,c_{1},c_{2}\in (-\infty,\,+\infty),$$

From this we obtain that the admissible function $\bar{x}(t)=t^{3},\,\, t\in [-1,\,1]$ is an extremal problem (5.1).

It is not difficult to assert that the admissible extremal
 $\bar{x}(t)=t^{3},\,\, t\in [-1,\,1]$ cannot be surrounded by a field, since the only one parametric family of extremals $x(t)=ct^{3},\,\, c\in (-\infty,\,+\infty), \, t\in [-1,\,1]$, containing $\bar{x}(t)=t^{3},\,\, t\in [-1,\,1]$ cannot cover the region inside which lies a point of the form $A(0,\,m)$, where $m\in(-\infty,\,+\infty)$.

It is obvious that the Jacobi-Weierstrass \cite{8} sufficient conditions can not be applied.

However, the application of Theorem 4.4 shows that the admissible extremal $\bar{x}(t)=t^{3},\,\,t\in[-1\,1]$ is an absolute minimum in problem (5.1).

\textbf{Example 5.2.} 
\begin{equation}
\int^{\frac{\pi}{2}}_{0}({\dot{x}}^{2}-x^{2})dt\rightarrow min,\,\, x(0)=1,\,x(\frac{\pi}{2})=0,
\end{equation} 
where $L(t,x,\dot{x})=\dot{x}^{2}-x^{2}$.

It is easy to see that the admissible function $\bar{x} (t)=\cos(t),\, t\in [0,\,\frac{\pi}{2}]$ is an extremal problem (5.3). Let us investigate the minimum of the admissible function $\bar{x} (t)=\cos(t),\, t\in [0,\,\frac{\pi}{2}]$.
)
We apply Theorem 4.2. The inequality (4.2) for the problem (5.3) along $\bar{x} (t)=\cos(t)$ have the form:

\begin{equation}
\xi^{2}+Q(t)\xi \eta+\frac{1}{2}(\dot{Q}(t)-2)\eta^{2}\geq 0,\,\,\forall \xi, \eta \in R, \,\,\,t\in[0,\,\frac{\pi}{2})
\end{equation} 

Solving the Ricati equation of the form

\begin{equation}
2\dot{Q}(t)=4+Q^{2}(t)
\end{equation} 
we obtain that the function $Q(t)=2\tan(t), \,\,t\in[0,\,\frac{\pi}{2})$ is a solution of equation (5.5). Further, it is clear that for $Q(t)=2\tan(t), \,\,t\in[0,\,\frac{\pi}{2})$ assumptions (3.2) are satisfied at the point $\frac{\pi}{2}$.

Therefore, the admissible extremal $ \bar{x}(t)= \cos(t),\,\, t\in[0,\,\frac{\pi}{2}]$ according to Theorem 4.2 is an absolute minimum in problem (5.3).

It is evident that the matrix function $Q(\cdot)$ (see Theorem 4.1 and Theorem 4.2) is generally defined as a solution to the matrix Riccata equation.

Note that in this paper the set $KC^{1}([t_{0},\,t_{1}],\,R^{n})$ is a class of admissible functions. This assumption in some cases limits the possibilities of using the sufficient minimum conditions obtained here. Therefore, it is natural to try to extend the class  $KC^{1}([t_{0},\,t_{1}],\,R^{n})$, for example, to the class of absolutely continuous functions.

To show the effectiveness of Theorem 4.1, it is enough to consider the following example.

\textbf{Example 5.3.}

\begin{align}
\int^{t_1}_{0}(x^{4/3}-x^{5/3}\dot{x}^{2})dt\rightarrow min,
\,\,\,\, x(0)=x(t_1)=0,
\end{align}
where $t_1>0$ and $L(t, x, \dot{x})=x^{4/3}-x^{5/3}\dot{x}^{2}$.
It's not difficult to find out that $\overline{x}(t)=0, t\in [0, t_1]$ is an admissible extremal for the problem (5.6) and along it the Weierstrass  conditions and the Legendre conditions degenerate, i.e. 
$$E(L)(t, \overline{x}(t), \dot{\overline{x}}(t), \dot{\overline{x}}(t)+\Delta \dot{x}(t))=0, \,\,L_{\dot{x}\dot{x}}(t, \overline{x}(t), \dot{\overline{x}}(t))=0, \,\,\, \forall t \in [0, t_1],$$ 
$$\forall \Delta x(t) \in KC^1 ([0, t_1], R).$$

Therefore, it is clear that the Weiestrass and Jacobi conditions [1-4] are not effective in the study of $\overline{x}(\cdot)$ in the problem (5.6). However, by applying Theorem 4.1 and taking $Q(\cdot)=0$ and $\delta = \frac{1}{2}$ we easily come to the conclusion that $\overline{x}(\cdot)=0$ is a weak local minimum in the problem (5.6).

We considered that obtaining analogs of the results of Section 4 for more general problems of the calculus of variations and for optimal control problems is promising. Undoubtedly, the possibilities of applying the concept of a set of integrands are extremely wide.


\end{document}